\magnification=\magstep1
\input amstex
\documentstyle{amsppt}
\redefine\R{\Bbb R}
\define\non{\text{\tt non}}
\define\cov{\text{\tt cov}}
\define\add{\text{\tt add}}
\define\cof{\text{\tt cof}}
\define\ctble{\text {\tt ctble}}
\define\dom{\text{\rm dom}}
\define\bdd{\text{\tt bounded}}
\redefine\S{\Bbb S}

\redefine\B{\text{\rm Borel}}

\define\rng{\text{\rm rng}}
\pageheight{7in}
\topmatter
\title
Duality chipped
\endtitle
\author
Jind\v rich Zapletal
\endauthor
\affil
University of Florida
\endaffil
\abstract
Whenever $I$ is a projectively generated projectively defined ideal, if ZFC+large cardinals $\vdash\cov(I)=\frak c$ then
ZFC+large cardinals $\vdash\non(I)<\aleph_4.$
\endabstract
\thanks
The author is partially supported by grants GA \v CR 201-00-1466 and NSF DMS-0071437.
\endthanks
\address
Department of Mathematics, University of Florida, Gainesville FL 32611
\endaddress
\email
zapletal\@math.ufl.edu
\endemail
\endtopmatter

\document

\head {0. Introduction}\endhead

The duality is a well-known phenomenon in the theory of cardinal invariants of the continuum \cite {B}. Given a $\sigma$-ideal $I$ on 
the reals, the invariant $\cov(I)$ is dual to $\non(I)$ and the invariant $\cof(I)$ is dual to $\add(I)$. The duality
heuristic states that ZFC proves an inequality of the form $x(I)\geq y(J)$ just in case it proves the inequality
$x^{dual}(I)\leq y^{dual}(J),$ where $x,y$ are among the symbols $\cov,\non,\cof,\add.$ While the proofs are seldom literal
translations of each other, still this heuristic is one of the most reliable experimentally established rules in the field.

The purpose of this paper is to prove

\proclaim {0.1. Theorem}
Whenever $I$ is a projectively generated projectively defined ideal, if ZFC+large cardinals $\vdash\cov(I)=\frak c$ then
ZFC+large cardinals $\vdash\non(I)<\aleph_4.$
\endproclaim

Here a $\sigma$-ideal $I$ is projectively generated if there is a number $n$ such that every set in $I$ is included in
a boldface $\Sigma^1_n$ set in $I$ and it is projectively defined if the set of all codes for boldface $\Sigma^1_n$ sets
in $I$ is itself boldface $\Sigma^1_m$ for some number $m$. In order for the theorem to make literal sense the latter set must be
actually lightface $\Sigma^1_m$ and the definition of this set is the definition of the ideal $I$. The reader is
invited to formulate and prove the boldface versions of the theorem.

To see how the theorem relates to the duality heuristic, note that $\frak c=\cof(\ctble)$ and $\aleph_1=\add(\ctble)$.
Thus the heuristic calls for replacing the $\aleph_4$ in the theorem by $\aleph_2$. However this is impossible:

\proclaim {0.2. Example}
There is a projectively generated projectively defined $\sigma$-ideal $I$ on the reals such that ZFC $\vdash\cov(I)=\frak c$
and ZFC+PFA $\vdash\non(I)=\aleph_2$.
\endproclaim

Theorem 0.1 continues to hold if the pair $\cov,\non$ is replaced by any other pair of dual symbols. Among the four resulting
statements one is true vacuously: if $I$ is a projectively generated projectively defined $\sigma$-ideal then ZFC+LC $\not\vdash\add(I)\geq
\cov$(meager) and so ZFC+LC $\not\vdash \add(I)=\frak c.$

The questions raised by Theorem 0.1 run in three different directions:

\definition {0.3. Question}
Can $\aleph_4$ in Theorem 0.1 be replaced by $\aleph_3$?
\enddefinition

As far as I can see this relates to the ZFC provability of partial squares on $\omega_2$ and therefore should be quite hard.

\definition {0.4. Question}
The duality heuristic asserts an equiprovability. Is there a version of Theorem 0.1 with the implication reversed?
\enddefinition

The converse of Theorem 0.1 is quite easily false. There should be a general theorem saying that the converse is true in many cases,
however my present understanding of the issues involved falls far short of proving such a result.

\definition {0.5. Question}
Are there similar theorems for invariants other than $\frak c$?
\enddefinition

The answer here is yes, but the theorems I can prove at this time have a certain provisional character to them. I will prove

\proclaim {0.6. Theorem}
Whenever $I$ is a projectively generated projectively defined ideal, if ZFC+large cardinals $\vdash\cov(I)\geq\frak d$ then
ZFC+large cardinals $\vdash\frak b\geq\aleph_3\to\non(I)\leq\frak b.$
\endproclaim

The notation used in this paper follows the set theoretic standard of \cite {J}. The large cardinal hypothesis used can be everywhere specified to
be ``there are $\omega_1$ Woodin cardinals'' by unpublished work of Neeman \cite {N}; its use below is denoted by a simple LC.
For an ordinal $\alpha$ and a set $C$ of ordinals the
symbol $\alpha^{+C}$ denotes the smallest ordinal in $C$ above $\alpha.$ For a $\sigma$-ideal $I$ the symbols $\cov(I),\non(I),\add(I)$ and $\cof(I)$
denote in turn the smallest number of $I$ small sets necessary to cover the real line, the smallest size of an $I$-positive set, the smallest possible
size of a collection $A\subset I$ with $\bigcup A\notin I$ and the smallest possible size of a collection $A\subset I$ cofinal
in the order $\langle I, \subset\rangle.$ The theorems are proved in the first section and Example 0.2
is proved in the second section.

\head {1. The theorems} \endhead

The argument of the main theorems contains two ingredients: the analysis of Sacks and Miller forcing from the determinacy point of view,
and a ZFC approaching sequences in regular cardinals.

\subhead{1.1. The Sacks forcing}\endsubhead

\definition {1.1.1. Definition}
Given a countable ordinal $\alpha$ let $\ctble^\alpha$ be the ideal of those sets $A\subset\R^\alpha$ for which player I
has a winning strategy in the game $G_A.$ The game has $\alpha$ many rounds and in the $\beta$-th round player I indicates a code for
a countable set of reals and player II indicates a real $r_\beta$ that does not belong to this countable set. Player II
wins if the sequence $\langle r_\beta:\beta\in\alpha\rangle$ belongs to the set $A.$
\enddefinition

\proclaim {1.1.2. Fact} (ZFC+LC)
For every countable ordinal $\alpha$ the countable support iteration of $\alpha$ many Sacks reals
is forcing equivalent to the poset $\B(\R^\alpha)\setminus\ctble^\alpha$ ordered by inclusion, which is
a dense subset of the poset of all projective subsets of $\R^\alpha$ minus the ideal $\ctble^\alpha.$
\endproclaim

This is the contents of \cite {Z, Section 1}.

\proclaim {1.1.3. Fact} (ZFC+LC)
For every countable ordinal $\alpha$ and every Borel positive set $A\subset\R^\alpha$ there is a Borel function
$g:\R^\alpha\to A$ such that the $g$-preimages of $\ctble^\alpha$-small sets are $\ctble^\alpha$-small.
\endproclaim

This is clear using the dense set $\S^\alpha$ of $\B(\R)\setminus\ctble^\alpha$ presented in \cite {Z, Definition 1.1.1}. It is convenient here
to put $\R=2^\omega$ and observe that every perfect set of reals there is in a canonical one-to-one correspondence with $\R$.
This observation naturally extends to yield a canonical level-preserving Borel bijection $g$ between $\R^\alpha$ and any Borel set $A\in\S^\alpha$,
meaning that $g(\vec r)\restriction\beta$ depends only on $\vec r\restriction\beta$ for every ordinal $\beta\in\alpha$.
For such a bijection the preimages of $\ctble^\alpha$-small sets are $\ctble^\alpha$-small. 

\proclaim {1.1.4. Fact} 
Let $I$ be a projectively generated projectively definable ideal. If ZFC+LC $\vdash\cov(I)=\frak c$ then ZFC+LC $\vdash\phi(I)$
where $\phi(I)=$``for some countable ordinal $\alpha$ and some Borel function $f:\R^\alpha\to\R$ the $f$-preimages of $I$-small
sets are $\ctble^\alpha$-small''.
\endproclaim

This fact is implicit in \cite {Z} and it is best proved in the contrapositive. Arguing in ZFC+LC, if in a countable support iteration
of Sacks forcing there appears a name for a real that falls out of every ground model coded set in the ideal $I,$ such a name
must appear at some countable stage $\alpha$ \cite {Z, Lemma 2.2.1}. Such a name takes form of a Borel function 
$f:B\to\R$ for some Borel $\ctble^\alpha$-positive set, for which the $f$-preimages of $I$-small
sets are $\ctble^\alpha$-small. Composing $f$ with the function $g$ from the previous Fact if necessary it is possible to arrange $\dom(f)=
\R^\alpha.$ Thus, if ZFC+LC $\not\vdash\phi(I)$ then in the consistent theory ZFC+LC+$\lnot\phi(I)$ it is possible to conclude that
in the iterated Sacks model the continuum is covered by the ground model coded sets in the ideal $I,$ meaning that the theory
ZFC+LC+$\cov(I)<\frak c$ is consistent.

\subhead {1.2. The combinatorics of regular cardinals}\endsubhead

The combinatorial part of the proof relates to the approachability and club guessing principles of Shelah.

\proclaim {1.2.1. Fact}
For every regular cardinal $\kappa\geq\aleph_3$ there is an approaching sequence: a sequence $\langle C_\delta:\delta\in\kappa\rangle$
together with a stationary set $S\subset\{\alpha\in\kappa:cf(\alpha)=\omega_1\}$ such that

\roster
\item $C_\alpha\subset\alpha$ is a closed set
\item for every ordinal $\beta\in C_\alpha$ except the maximum there is $\gamma\in\beta^{+C_\alpha}$ such that $C_\gamma=C_\alpha\cap\beta+1$
\item for every ordinal $\alpha\in S$ $C_\alpha\subset\alpha$ is unbounded of ordertype $\omega_1.$
\endroster
\endproclaim

I could not find a clean proof of this relatively well-known fact in print. It can be assembled from \cite {S III, Claim 2.14}
and \cite {M, Lemma 2.1.1, 2.1.2}. It is not known whether approaching sequences exist
at $\kappa=\aleph_2$ in ZFC, and this is quite annoying for the application in this paper. Clearly $\square_{\aleph_2}$ implies the
existence of an approaching sequence for $\aleph_2.$

Using the standard transfinite adjustment process of \cite {S} it is possible to thin out any approaching sequence to one with a club
guessing property:

\proclaim {1.2.2. Fact}
For every regular cardinal $\kappa\geq\aleph_3$ there is a club guessing approaching sequence: a sequence $\langle C_\delta:\delta\in\kappa\rangle$
together with a stationary set $S\subset\{\alpha\in\kappa:cf(\alpha)=\omega_1\}$ such that (1-3) of the previous Fact hold, together
with

\roster
\item"{(4)}" for every closed unbounded set $E\subset\kappa$ there is $\delta\in S$ such that $C_\delta\subset E.$
\endroster
\endproclaim

\subhead {1.3. The proofs of Theorems 0.1 and 0.6}\endsubhead

The key claim connecting the previous two subsections is

\proclaim {1.3.1. Lemma}
For every countable ordinal $\alpha,$ $\non(\ctble^\alpha)<\aleph_4$.
\endproclaim

\demo {Proof}
Fix an arbitrary countable ordinal $\alpha.$ 
If $\frak c<\aleph_3$ then we are done since $\frak c=|\R^\alpha|<\aleph_4$ and $\R^\alpha$ is certainly a $\ctble^\alpha$-positive
set. If $\frak c\geq\aleph_3$ fix a sequence $\langle r_\gamma:\gamma\in\omega_3\rangle$ of distinct reals and a club guessing
approaching sequence $\langle C_\delta:
\delta\in\omega_3\rangle$ together with the relevant stationary set $S\subset\{\delta\in\omega_3:cf(\delta)=\omega_1\}$. 
For each ordinal $\delta\in\omega_3$ let $\vec r_\delta$
be the sequence of reals $\langle r_\gamma:\gamma$ a nonaccumulation point of the set $C_\delta\rangle$ or the first $\alpha$ elements of this sequence,
whichever is shorter. These sequences will be viewed alternately as indexed by ordinals below $\alpha$ or by the nonaccumulation points
in the set $C_\delta$, whichever is more convenient. The proof will be complete once I show that the set
$A=\{\vec r_\delta:\delta\in S\}$ is $\ctble^\alpha$-positive.

Well, let $\sigma$ be a strategy in the game $G_A$ for player I. 
An ordinal $\delta\in S$ must be found such that the sequence $\vec r_\delta$ is a legal counterplay against
the strategy $\sigma$. First find a continuous tower $T$ of height $\omega_3$ of submodels of some large structure containing 
all the objects mentioned so far, such that every model $M\in T$ has size $\aleph_2$ and $M\cap\omega_3\in\omega_3.$
There is an ordinal $\delta\in S$ such that $C_\delta$ is a subset of the club $\gamma\in\omega_3:\exists M\in T\ \gamma=M\cap\omega_3\}.$

Now the sequence $\vec r_\delta$ is the desired legal counterplay against the strategy $\sigma.$ To see this, choose an arbitrary 
ordinal $\beta\in\dom(\vec r)$; it must be the case that the real $r_\beta$ does not belong to the countable set $\sigma(\vec r_\delta\restriction\beta)$
that the strategy $\sigma$ produces after player II has played the reals on the sequence $\vec r_\delta\restriction\beta$ in their turn.
There is a model $M\in T$ such that $\beta=M\cap\omega_3.$ For some ordinal $\gamma\in\beta,$ $C_\delta\cap\beta=C_\gamma$ and so $\vec r_\gamma=
\vec r_\delta\restriction\beta.$ All the objects $\gamma, \vec r_\gamma$ and $\sigma(\vec r_\gamma)$ are in the model $M$ and since 
$\sigma(\vec r_\gamma)$ is a countable set of reals, $\sigma(\vec r_\gamma)=\sigma(\vec r_\delta\restriction\beta)\subset M.$
However, $r_\beta\notin M$ and so $r_\beta\notin\sigma(\vec r_\delta\restriction\beta)$ as desired. The lemma follows.
\qed
\enddemo

Theorem 0.1 now quickly follows. Suppose that ZFC+LC $\vdash\cov(I)=\frak c.$ By Fact 1.1.4, ZFC+LC $\vdash$ for some countable ordinal
$\alpha$ and a Borel function $f:\R^\alpha\to\R$, the $f$-preimages of $I$-small sets are $\ctble^\alpha$-small.
Now argue in ZFC+LC. Choose an ordinal $\alpha$ and a function $f$ as above and choose a $\ctble^\alpha$-positive
set $X\subset\R^\alpha$ of size $<\aleph_4$ by Lemma 1.3.1. Then the set $f''X$ must be $I$-positive of size $<\aleph_4$, and
so $\non(I)<\aleph_4.$

The proof of Theorem 0.6 is essentially the same. I will indicate only the main changes.

\definition {1.3.2. Definition}
Given a countable ordinal $\alpha$ let $\bdd^\alpha$ be the ideal of those sets $A\subset(\omega^\omega)^\alpha$ for which player I
has a winning strategy in the game $H_A.$ The game has $\alpha$ many rounds and in the $\beta$-th round player I indicates a function
$k\in\omega^\omega$ and player II indicates a function $h_\beta$ that is not modulo finite dominated by $k$. Player II
wins if the sequence $\langle h_\beta:\beta\in\alpha\rangle$ belongs to the set $A.$
\enddefinition

The countable support iteration  of countable length of the Miller forcing is forcing-equivalent to the poset $\B((\omega^\omega)^\alpha)\setminus
\bdd^\alpha$ as proved in \cite {Z}. The proof of Theorem 0.6 now proceeds just as above, the main point being

\proclaim {1.3.3. Lemma}
If $\frak b\geq\aleph_3$ then $\non(\bdd^\alpha)=\frak b$ for all countable ordinals $\alpha.$
\endproclaim

Note that $\frak b=\non(\bdd^1)$ by the definitions, and as the ordinal $\alpha$ increases, so could the invariant $\non(\bdd^\alpha)$.
The Lemma shows that this does not happen as long as $\frak b$ is large enough. The proof is almost the same as the argument
for Lemma 1.3.1 with $\aleph_3$ replaced by $\frak b$. Note that $\frak b$ is necessarily regular and so the Fact 1.2.1 is applicable.
In the beginning choose a modulo finite increasing unbounded sequence $\langle h_\delta:\delta\in\frak b\rangle$ of functions in $\omega^\omega$
 in place
of the distinct reals $r_\delta$ and proceed as before.

\head {2. The example}\endhead

Let $F:\R^2\to\R$ be a function. Define a $\sigma$-ideal $I_F$ on $\R^\omega$ as the one generated by sets $A_X\subset\R^\omega,$
where $X$ ranges over all countable sets of reals and $A_X=\{\vec r\in\R^\omega:$ for some number $n$ the real
$\vec r(n)$ belongs to the $F$-closure of the set $X\cup\{\vec r(m):m\in n\}\}.$

Clearly, $\cov(I_F)=\frak c.$ For if $\{X_\alpha:\alpha\in\kappa\}$ is a collection of $\kappa<\frak c$ many countable sets
of reals, a sequence $\vec r\in\R^\omega\setminus\bigcup_{\alpha\in\kappa}A_{X_\alpha}$ can be constructed by induction.
Note that given the reals $\vec r(m):m\in n,$ the $F$-closure of the set $\bigcup_{\alpha\in\kappa}X_\alpha\cup\{\vec r(m):m\in n\}$
has size $\leq\kappa<\frak c$ and thus it is possible to choose a real $\vec r(n)$ so that it does not belong to this set. 
If this is done at each $n\in\omega,$ the resulting sequence $\vec r$ will fall outside of all sets $A_{X_\alpha},\alpha\in\kappa.$
The argument for Example 0.2 will be complete once I find a (lightface) Borel function $F:\R^2\to\R$ such that ZFC+PFA
$\vdash\non(I_F)=\aleph_2.$

For the remainder of this section put $\R=[\omega]^{\aleph_0}.$ Choose a Borel injection $g:\R\to\R$ whose range consists of
mutually almost disjoint sets of reals. Choose a Borel function $h:\R\to\R^\omega$ whose range consists of injections
with mutually disjoint ranges--$r\neq s\to \rng(h(s))\cap\rng(h(r))=0$. 
This is a setup for an almost disjoint coding: every real $x$ defines a function $F_x:
\R\to\R$ by $n\in F_x(y)\leftrightarrow x\cap g(h(y)(n))$ is infinite. Let $F(x,y)=F_x(y)$.

In order to show that PFA$\to\non(I_F)=\aleph_2$ I will produce a proper forcing $P\Vdash(\R^\omega)^V\in I_F.$ Assume
without loss of generality that the continuum hypothesis holds or else force it by a $\sigma$-closed notion of forcing.
Let $\langle r_\alpha:\alpha\in\omega_1\rangle$ be a one-to-one enumeration of the reals. 
The poset $P$ will be an iteration $P_0*\dot P_1$ where

\roster
\item $P_0$ is the Baumgartner's forcing for adding a closed unbounded subset $C$ of $\omega_1$ with finite conditions
\item $P_1$ is a variation of Solovay's almost disjoint c.c.c. coding adding reals $x_n:n\in\omega$ such that for every
ordinal $\alpha\in\omega_1,$ $\{r_\beta:\beta\in\alpha^{+C}\}\subset F''\{\langle x_n,r_\alpha\rangle:n\in\omega\}.$
\endroster

In the end, let $X=\{x_n:n\in\omega\}$ and observe that $P\Vdash(\R^\omega)^V\subset A_X.$ For whenever $\langle r_{\alpha_n}:
n\in\omega\rangle$ is a ground model $\omega$-sequence of reals then by an elementary genericity argument for $P_0$ there must 
be distinct integers $m\in n$ such that the ordinals $\alpha_n, \alpha_m$ belong to the same hole in the set $C\subset\omega_1.$
But then the real $r_{\alpha_n}$ is in the $F$ closure of the set $X\cup\{r_m\}$ and $\langle r_{\alpha_n}:
n\in\omega\rangle\in A_X$ as desired.

To complete the proof of Theorem 0.2 I just have to demonstrate how the forcing $P_1$ is obtained. A condition in $P_1$ is
a finite sequence $\langle f, a_0, a_1, \dots a_m\rangle$ where $f:\omega_1\times\omega\to\omega_1$ is a finite function
such that for every pair $(\alpha, n)$ in its domain it is the case that $n\leq m$ and $f(\alpha, n)\in\alpha^{+C}$, and $a_n:n\leq m$ are finite
subsets of $\omega.$ The information carried by such a condition is that each $a_n$ is an initial segment of the future set $\dot x_n$,
for every pair $(\alpha, n)$ in the domain of $f$ the value $F(\dot x_n,r_\alpha)$ will be $r_{f(\alpha, n)}$ and to secure this fact,
for each number $k\in m\setminus r_{f(\alpha, n)}$ the finite set $\dot x_n\cap g(h(r_\alpha)(k))$ will be a subset of 
$a_n$. The ordering is defined accordingly:
$\langle g, b_0, b_1, \dots b_l\rangle\leq \langle f, a_0, a_1, \subset a_m\rangle$ if $f\subset g,$ $m\leq l$, for each number $n\leq m$
the set $a_n$ is an initial segment of $b_m$ and whenever $(\alpha, n)\in\dom(f)$ and $k\in m\setminus r_{f(\alpha, n)}$ then
$g(h(r_\alpha)(k))\cap a_n=g(h(r_\alpha)(k))\cap b_n$. Standard density arguments will show that indeed, putting $\dot x_n$
to be the name for the union of all sets $a_n$ in the conditions in the generic filter, $P_1\Vdash$ the sets $\{x_n:n\in\omega\}$
have the property from (2) above. The only thing left to verify is the c.c.c. of the poset $P_1$ and that is completely standard.
Using a $\Delta$-system argument, in any uncountable subset of $P_1$ it is possible to find conditions $\langle f, a_0, \dots a_m\rangle,
\langle g, a_0, \dots a_m\rangle$ such that the finite subsets of $\omega$ mentioned in the conditions will be the same, and
$f\cup g$ will be a function. Clearly such two conditions are compatible and their lower bound is  $\langle f\cup g, a_0, \dots a_m\rangle$.

As the last remark, one can show in ZFC that $\non(I_F)\leq\non\ctble^\omega<\aleph_3$ using a weak club-guessing principle on $\aleph_2$
provable in ZFC. However, the argument is not available for the versions of the ideal $I_F$ for sequences of reals of length $\omega^2$ or more.

\enddocument